\documentclass{amsart}
\newtheorem{Theorem}{Theorem}[section]
\newtheorem{Definition}[Theorem]{Definition}
\newtheorem{Proposition}[Theorem]{Proposition}

\newtheorem{Lemma}[Theorem]{Lemma}
\newtheorem{Corollary}[Theorem]{Corollary}
\theoremstyle{remark}
\newtheorem{Example}[Theorem]{Example}

\def\il{\int\limits_}

\def\eps{\varepsilon}

\def\om{\omega}

\def\bd{\partial}

\def\sbs{\subset}

\def\supp{\operatorname{supp}}

\def\be{\begin{enumerate}}
\def\ee{\end{enumerate}}
\def\bT{\begin{Theorem}}
\def\eT{\end{Theorem}}
\def\bP{\begin{Proposition}}
\def\eP{\end{Proposition}}
\def\bD{\begin{Definition}}
\def\eD{\end{Definition}}
\def\bE{\begin{Example}}
\def\eE{\end{Example}}
\def\bL{\begin{Lemma}}
\def\eL{\end{Lemma}}
\def\bC{\begin{Corollary}}
\def\eC{\end{Corollary}}

\def\M{{\mathcal M}}
\def\J{{\mathcal J}}
\def\rO{{\mathcal O}}

\begin{document}
\title{The Dirichlet Problem for Harmonic Functions on Compact Sets}
\author{Tony L. Perkins}

\address{ Department of Mathematics,  215 Carnegie Building,
Syracuse University,  Syracuse, NY 13244-1150}
\begin{abstract}
The primary goal of this paper is to study the Dirichlet problem on a compact set $K\subset \mathbb{R}^n$. Initially we consider the space $H(K)$ of functions on $K$ which can be uniformly approximated by functions harmonic in a neighborhood of $K$ as possible solutions.  As in the classical theory, our Theorem~\ref{T:regDP} shows $C(\partial_fK)\cong H(K)$ for compact sets with $\partial_fK$ closed, where $\partial_fK$ is the fine boundary of $K$. However, in general a continuous solution cannot be expected even for
continuous data on $\partial_fK$ as illustrated by Theorem~\ref{T:regDP}.
Consequently, we show that the solution can be found in a class of
finely harmonic functions. Moreover
by Theorem~\ref{T:qdp}, in complete analogy with the classical
situation, this class is isometrically isomorphic to
$C_b(\partial_fK)$ for all compact sets $K$.
\end{abstract}
\date{\today}
\maketitle
\section{Introduction}
The Dirichlet problem for harmonic functions on domains in
$\mathbb{R}^n$ is not only important by itself but by its
influence on potential theory.  Many now standard notions, e.g.
regular points, fine topology, etc., first appeared in the
study of this problem.   The main goal of the present paper is
to extend the classic theory to compact sets $K\subset
\mathbb{R}^n$.

One possible extension can be found in the abstract theory of balayage spaces, see \cite{BH86,H85}.  However we feel that the gain in transparency following from a direct geometric approach more than justifies the use of new techniques.

The Dirichlet problem can be thought of as having two
components; the data set and the data itself.  One uses an
initial function defined on the data set to construct a
solution (a harmonic function) on the rest of the domain which
must have a prescribed regularity as it approaches the data
set. Classically, the data set is taken to be the topological
boundary of the domain.  One of the main goals of this paper is
to establish that the natural choice for the data set on
compact sets is the fine boundary of $K$, $\partial_fK$, which is shown by Lemma \ref{L:Pfine} to be the Choquet boundary of $K$ with respect to subharmonic functions on $K$. We limit
ourselves to initial functions that are continuous and bounded on $\partial_fK$ as
in the classical case.

In Section~\ref{S:HSub}, we introduce Jensen measures as our main tool and begin extending potential theory to compact sets $K\subset \mathbb{R}^n$ by defining harmonic functions and subharmonic functions on $K$.  We devote Section~\ref{S:Hmeas} to the construction and study of harmonic measure on compact sets. The harmonic measure on $K$ is shown to be a maximal Jensen measure.  This is used to see the important fact (Corollary \ref{C:conc}) that harmonic measures are concentrated on the fine boundary. In Section~\ref{S:DP} we study the Dirichlet problem for compact sets.  As in the classical theory, our Theorem~\ref{T:regDP} shows $C(\partial_fK)\cong H(K)$ for a class of compact sets. However, in
general a continuous solution cannot be expected even for
continuous data on $\partial_fK$ as illustrated by Example~\ref{E:SSC}.
Consequently, we show that the solution can be found in the class of
finely harmonic functions introduced in this section. Moreover
by Theorem~\ref{T:qdp}, in complete analogy with the classical
situation, this class is isometrically isomorphic to
$C_b(\partial_fK)$ for all compact sets $K$.

\thanks{
It is a pleasure to thank Stephen J. Gardiner, Leonid Kovalev, and Gregory Verchota for stimulating discussions related to the topic of this article.  We are especially grateful to Eugene Poletsky for his excellent guidance and support.
}

\section{Basic Facts}
First some notation.  Let $\mathcal{M}(\Omega)$ denote the space of finite signed Radon measures on $\Omega\subset \mathbb{R}^n$ and $C_0(\mathbb{R}^n)$ will be the space of continuous functions on $\mathbb{R}^n$ which vanish at infinity.  We will often use $\mu(f)$ to denote $\int f \, d\mu$.
\subsection{Classical Potential Theory}
Let $D$ be an open set in $\mathbb{R}^n$, $n\ge 2$.  For any
$f\in C(\partial D)$, the \emph{Dirichlet problem} on $D$ is to
find a unique function $h$ which is harmonic on $D$ and
continuous on $\overline{D}$ such that $h|_{\partial D} = f$.
The function $f$ is commonly referred to as the \emph{boundary
data}, and the corresponding $h$ is said to be the
\emph{solution} of the Dirichlet problem on $D$ with boundary
data $f$.  The punctured disk in $\mathbb{R}^2$ is a fundamental example which
shows that the Dirichlet problem can not be solved for any
continuous boundary data.  However for a bounded open set $U$
the method of Perron allows one to assign a function which is
harmonic on $U$ to any continuous (or simply measurable)
boundary data.  Later the concept of a regular domain was
developed to establish the continuity of the Perron solution to
the boundary.  A bounded connected open set $D\subset \mathbb{R}^n$ is
a \emph{regular} domain if the Dirichlet problem is solvable on
$D$ for any continuous boundary data.  Therefore on a regular
domain, $C(\partial D)$ is isometrically isomorphic to $H(D)$,
the space of continuous functions on $\overline{D}$ which are
harmonic on $D$.  For any $f\in C(\partial D)$ let $h_f\in
H(D)$ denote the solution of the Dirichlet problem on $D$ with
boundary data $f$.  Let $z\in D$.  The point evaluation $H_z\colon f\mapsto h_f(z)$ is a
positive bounded linear functional on $C(\partial D)$.  By the
Riesz Representation Theorem, there is a Radon measure
$\om_D(z, \cdot)$ on $\partial D$ which represents $H_z$, that is
\[h_f(z) = \il{\partial D} f(\zeta)\; d\om_D(z, \zeta),\]
 for all $f\in C(\partial D)$. The measure $\om_D(z, \cdot)$ is
called the \emph{harmonic measure} of $D$ with barycenter at
$z$.  See \cite{AG01} for more
details on potential theory.

\subsection{Jensen Measures}
If $D$ is an open set in $\mathbb{R}^n$, we say that $\mu$ is a {\it Jensen measure} on $D$ with barycenter $z\in D$ if $\mu$
is a probability measure (a positive Radon measure of
unit mass) whose support is compactly contained in $D$ and for
every subharmonic function $f$ on $D$ the {\it sub-averaging
inequality} $f(z)\le\mu(f)$ holds.  The set of Jensen measures
on $D$ with barycenter $z\in D$ will be denoted $\J_z(D)$.

One could define the set of Jensen measures $\J_z^c(D)$ with respect to the continuous subharmonic functions on $D$.  However the following theorem shows that the set of Jensen measures would not be changed.

\bT\label{T:JensEQ} Let $D$ be a bounded open subset of $\mathbb{R}^n$.
For every $z\in D$, the sets $\J_z(D)$ and $\J^c_z(D)$ are equal.
\eT
\begin{proof}
Since it is clear that $\J_z(D) \subseteq \J^c_z(D)$ for all
$z\in D$, we will now show the reverse inclusion.

Pick some $z_0\in D$ and let $\mu\in \J^c_{z_0}(D)$.  Then we must show
$f(z_0)\le \mu(f)$ for every function $f$ which is subharmonic
on $D$.  The support of $\mu$ is compactly contained in $D$.

Since $f$ is subharmonic on $D$ we can find an decreasing sequence $\{f_n\}$ of continuous subharmonic functions which converge to $f$.  As $\mu\in \J_{z_0}^c(D)$ we have $f(z_0) \le \mu(f_n)$ for every $f_n$.   By the Lebesgue Monotone Convergence Theorem
it follows that $f(z_0)\le \mu(f)$.  Thus $\mu\in \J_{z_0}(D)$.
\end{proof}

Since $\J_z(D)=\J^c_z(D)$ for all $z\in D$, to check that
$\mu\in \J_z(D)$, it suffices to check that $\mu$ has the
sub-averaging property for every continuous subharmonic
function.

Examples of Jensen measures with barycenter at $z\in D$ include the Dirac measure at $z$, i.e. $\delta_z$, the
harmonic measure with barycenter at $z$ for any regular domain
which is compactly contained in $D$, and the average over any
ball (or sphere) centered at $z$ which is contained in $D$.
The following proposition of Cole and
Ransford~\cite[Proposition 2.1]{CR01} will demonstrate some
basic properties of sets of Jensen measures.

\bP\label{P:basic} Let $D_1$ and $D_2$ be open subsets of $\mathbb{R}^n$ with $D_1\subset D_2$.  Let $z\in D_1$.
\begin{enumerate}
\item[{\rm i.}] If $\mu\in \J_z(D_1)$ then also $\mu\in
    \J_z(D_2)$.
\item[{\rm ii.}] If $\mu\in \J_z(D_2)$ and
    $\supp(\mu)\subset D_1$, and if each bounded component
    of $\mathbb{R}^n\setminus D_1$ meets
    $\mathbb{R}^n\setminus D_2$, then $\mu\in \J_z(D_2)$.
\end{enumerate}
\eP

Jensen measures and subharmonic functions are, in a sense, dual to each other.  This duality is illustrated by the following
theorem of Cole and Ransford~\cite[Corollary 1.7]{CR97}.
\bT\label{T:duality}
Let $D$ be an open subset of $\mathbb{R}^n$ which possesses a
Green's function.  Let $\phi\colon D\rightarrow [-\infty,
\infty)$ be a Borel measurable function which is locally
bounded above.  Then, for each $z\in D$,
\[\sup\left\{v(z)\colon v\in S(D), v\le \phi\right\}=\inf\left\{\mu(\phi)\colon \mu \in \J_z(D)\right\},\]
where $S(D)$ denotes the set of subharmonic functions on $D$.
\eT

\subsection{Fine Topology}
The two books \cite{B71, F72} are classical references on
the fine topology and many books on potential theory contain
chapters on the topic, e.g. \cite[Chapter 7]{AG01}.

The \emph{fine topology} on $\mathbb{R}^n$ is the coarsest
topology on $\mathbb{R}^n$ such that all subharmonic functions
are continuous in the extended sense of functions taking values
in $[-\infty, \infty]$.

When referring to a topological concept in the fine topology we
will follow the standard policy of either using the words
``fine"  or  ``finely" prior to the topological concept or
attaching the letter $f$ to the associated symbol.
For example, the fine boundary of $K$, $\partial_fK$,
is the boundary of $K$ in the fine topology.  The fine topology
is strictly finer than the Euclidean topology.

Many of the key concepts of classical potential theory have
analogous definitions in relation to the fine topology.
Presently we will recall a few of them.  Relative to a finely
open set $V$ in $\mathbb{R}^n$ the \emph{harmonic measure}
$\delta_x^{V^c}$ is defined as the swept-out of the Dirac
measure $\delta_x$ on the complement of $V$.  A function $u$ is
said to be \emph{finely hyperharmonic} on a finely open set $U$
if it is lower finite, finely lower semicontinuous, and
\[ -\infty< \delta_x^{V^c}(u)\le u(x),\]
 for all $x\in V$ and all relatively compact finely open sets $V$ with fine closure contained in $U$.  A function $h$ is said to be \emph{finely harmonic} if $h$ and $-h$ are finely hyperharmonic.  Furthermore, the \emph{fine Dirichlet problem on $U$} for a finely continuous function $f$ defined on the fine boundary of a bounded finely open set $U$ consists of finding a finely harmonic extension of $f$ to $U$.  The development of the fine Dirichlet problem is quite similar to that of the classical.  In the seventies Fuglede~\cite{F72} establishes a Perron solution for the fine Dirichlet problem.  His~\cite[Theorem 14.6]{F72} shows that there exists a Perron solution $H_f^U$ which is finely harmonic on $U$ for any numerical function $f$ on $\partial_fU$ which is $\delta_x^{\partial_fU}$ integrable for every $x\in U$.  Furthermore~\cite[Theorem 14.6]{F72} provides us with the desired continuity at the boundary, i.e. that the fine limit of $H_f^U(x)$ tends to $f(y)$ as $x\in U$ goes to $y$ for every finely ``regular" boundary point $y\in \partial_fU$ at which $f$ is finely continuous.

\section{Harmonic and Subharmonic Functions on Compact Sets}\label{S:HSub}

We now begin our study of potential theory on compact sets.  For compact sets which are not connected, the Hausdorff property will allow us to reduce Dirichlet type problems on the compact set to solving separate problems on each connected component.  Therefore in what follows we will work on compact sets $K$ in $\mathbb{R}^n$ which need not be connected, with the understanding that we can always separate the problem by working the connected components of $K$ individually.

There are currently three equivalent ways to define harmonic and subharmonic functions on compact sets.
\bD[Exterior]
Let $H(K)$ (or $S(K)$) be the unform closures of  all functions in $C(K)$ which are restrictions of harmonic (resp. subharmonic) functions on a neighborhood of $K$.
\eD

\bD[Interior] One can define $H(K)$ and $S(K)$ as the subspaces of $C(K)$ consisting of functions which are finely harmonic (resp. finely superharmonic) on the fine interior of $K$.
\eD
The equivalence of these definitions of $H(K)$ was shown in \cite{DG74} and of $S(K)$ in \cite{BH75, BH78}.

For the third definition of $H(K)$ we must to extend the notion of Jensen measures to compact sets.

\bD We define the set
of Jensen measures on $K$ with barycenter at $z\in K$ as the
intersection of all the sets $\J_z(U)$, that
is
\[\J_z(K) = \bigcap_{K\subset U} \J_z(U),\]
where $U$ is any open set containing $K$.
\eD

Another definition of $H(K)$ was introduced in \cite{P97} using the notion of Jensen measures.

\bD[Via Jensen measures] The set $H(K)$ is the subspace of $C(K)$ consisting of functions $h$ such that $h(x)=\mu(h)$ for all $\mu\in\J_x(K)$ and $x\in K$.
\eD

It was shown in \cite{P97} that this definition is equivalent to the exterior definition above.

Our first lemma shows that this last construction of Poletsky extends to subharmonic functions in the ideal way.

\bL A function is in $S(K)$ if and only if it is continuous and
satisfies the subaveraging property with respect to every Jensen
measure on $K$, that is\[S(K) = \left\{f\in C(K) \colon f(z)\le \mu(f), \text{ for all } \mu \in \J_z(K) \text{ and every }  z\in K \right\}.\]
\eL
\begin{proof}
We use the exterior definition of $S(K)$ to show ``$\subseteq$".  Let $\{f_j\}$ be a sequence of subharmonic functions defined in a neighborhood of $K$ such that $\{f_j\}$ is converging uniformly to $f$.  Then $f_j(z)\le \mu(f_j)$ for any $\mu\in \J_z(K)$.  Since the convergence is uniform we have $f(z) \le \mu(f)$.

Now suppose that $f$ is in the set on the right.  The subaveraging condition implies that $f$ is finely subharmonic, and by assumption $f$ is continuous.  Therefore $f$ satisfies the interior definition of $S(K)$.
\end{proof}

Recall the (exterior) definition of $S(K)$ as the uniform limits of continuous functions subharmonic in neighborhoods of $K$.  The following proposition shows that the defining sequence for any function in $S(K)$ may be taken to be increasing. This result is a simple consequence of a duality theorem of Edwards.
\bP\label{P:SbhExt}
Every function in $S(K)$ is the limit of an
increasing sequence of continuous subharmonic functions defined
on neighborhoods of $K$.
\eP
\begin{proof}
Recall (see \cite[Theorem 1.2]{G78}and \cite{CR97}) Edwards Theorem states: If $p$ is a continuous function on $K$,
then for all $z\in K$ we have
\[Ep(z):=\sup\{f(z)\colon f\in S(K),\; f\le p\}=\inf\{\mu(p)\colon \mu\in \J_z(K)\}.\]
From the proof of this theorem it follows that $Ep$ is lower semicontinuous and is the limit of an increasing sequence of continuous subharmonic functions on neighborhoods of $K$.  The result follows by observing that $p=Ep$ whenever $p\in S(K)$.
\end{proof}

\section{Harmonic Measure on a Compact Set}\label{S:Hmeas}

To use the exterior definition of $H(K)$ we will commonly want to approximate $K$ by a decreasing sequence of regular domains. A decreasing sequence of regular domains $\{U_j\}$ is said to be \emph{converging} to $K$ if for every $\eps>0$ there is a $j_0$ such that $U_j$ lies in the $\eps$-neighborhood $K_\eps$ of $K$ when $j\ge j_0$.   Furthermore, we require that $U_{j+1}$ is
compactly contained in $U_j$, i.e. $\overline{U}_{j+1} \subset
U_j$, for all $j$.  The existence of such a sequence is provided by \cite[Prop 7.1]{H62}.

The next theorem will allow us to define a harmonic measure on
$K$. For a decreasing sequence of regular domains $\{U_j\}$, we will let $\om_{U_j}(z,\cdot)$ denote the harmonic measure on $U_j$ with barycenter at $z\in U_j$.

\bT\label{T:1}
If $\{U_j\}$ is a sequence of regular domains converging to a
compact set $K\sbs\mathbb{R}^n$, then for every $z\in K$ the
harmonic measures $\omega_{U_j}(z,\cdot)$ converge weak$^\ast$.
Furthermore, this limit does not depend on the choice of the
sequence of domains $\{U_j\}$.
\eT
\begin{proof} Since $\om_{U_j}$ are measures of unit mass supported
on a compact set in $\mathbb{R}^n$, by Alaoglu's Theorem they
must have a limit point. To show that this point is unique it
suffices to show that for every $z\in K$ the limit
\begin{equation}\label{e:wchm}
\lim_{j\to\infty}\il{\partial U_j} u
(\zeta) \; d\omega_{U_j}(z, \zeta)
\end{equation}
exists for every $u\in C(\overline{U}_1)$.

First, we show the limit in (\ref{e:wchm}) exists when $u$ is
continuous and subharmonic in a neighborhood of $K$. The
solution $u_j$ of the Dirichlet problem on $U_j$ with boundary
value $u$ is equal to
\[ u_j(z) = \il{\partial U_j} u (\zeta) \;d\omega_{U_j}(z, \zeta).\]
Since $u$ is subharmonic, we have $u_j \ge u$ on $U_j$.  Then as
$u_{j+1}=u $ on $\partial U_{j+1}$ and $u_j \ge u = u_{j+1}$ on
$\partial U_{j+1}$, the maximum principle for harmonic
functions implies that $u_j  \ge u_{j+1} $ on $U_{j+1}$.  Thus
$\{ u_j \}$ is a decreasing sequence  on $K$ and we see that
for every $z\in K$ the limit in (\ref{e:wchm}) exists.

If $u\in C^2(\overline{U}_1)$, then we may represent $u$ as a
difference of two $C^2(\overline{U}_1)$ functions which are
subharmonic on $U_1$. By the argument
above the limit in (\ref{e:wchm}) exists.

Since $C^2(\overline{U}_1)$ is dense in $C(\overline{U}_1)$ we
see that the limit in (\ref{e:wchm}) always exists.
\end{proof}

\bD
We define the harmonic measure $\om_K(z,\cdot)$ on a compact set $K$ with $z\in K$ as the weak$^\ast$ limit of $\omega_{U_j}(z,\cdot)$ chosen as above.
\eD

To use this definition for the Dirichlet problem we must check that the support of $\om_K(z,\cdot)$ lies on the boundary of $K$.  Actually in Section \ref{S:Boundary} we will be able to give more specific information about $\om_K(z, \cdot)$, see Corollary \ref{C:conc}.

\bL
The support of $\om_K(z, \cdot)$ is contained in $\partial K$.
\eL
\begin{proof}
Let $W$ be a neighborhood of $\partial K$. Let $\{U_j\}$ be a
sequence of domains converging to $K$ and take a sequence
$z_j\in\partial U_j$.  Then there exists a subsequence
$\{z_{j_k}\}$ which must be converging to some $z_0\in K$. As
$z_j\in\partial U_j$, then $z_j$ is not in $K$. Therefore the
limit of $z_{j_k}$ cannot be in the interior of $K$.  Thus
$z_0$ is in $\partial K\subset W$. Consequently, there is $j_0$
such that $\partial U_j\subset W$ for each $j\ge j_0$,

Let $x\in \mathbb{R}^n\setminus \partial K$ and take $W$ to be
a neighborhood of $\partial K$ so that $x\not\in \overline{W}$.
There is an $r>0$ so that $\overline{B(x,r)}\cap \overline{W} = \emptyset.$  Since $\omega_{U_j}(z, \cdot)$ has support on $\partial U_j$, which is contained is $\overline{W}$ for large $j$, we have $\omega_{U_j}(z,B(x,r))=0$. Since $B(x, r)$ is open, the Portmanteau Theorem shows
\[\liminf_{j\rightarrow \infty}\omega_{U_j}(z,B(x,r)) \ge \omega_K(z,B(x,r)). \]
Hence $\om_K(z,B(x,r))=0$ and $x$ is not in the support of $\om_K(z, \cdot)$.
\end{proof}

The following theorem brings our study back to the topic of Jensen measures.

\bT\label{T:HJens}
The harmonic measure on $K$ is a Jensen measure on $K$.
\eT
\begin{proof}
Since $\om_K(z, \cdot)$ is defined as the weak$^\ast$ limit of
probability measures, $\om_K(z, \cdot)$ is a probability
measure.

Recall that for $z\in K$ we have defined
$\J_z(K)=\cap\J_z(U)$, where $K\subset U$.  However it is sufficient to see that $\J_z(K) = \cap\J_z(U_J)$ where $\{U_j\}$ is any sequence of
domains converging to $K$.  We will show $\om_K(z,
\cdot)\in\J_z(U_j)$ for all $j$.

Pick some $j$. Then let $f$ be a continuous subharmonic
function on $U_j$.  Then \[f(z)\le \il{\partial U_l} f(\zeta)\;
d\om_{U_l}(z, \zeta),\] for all $l > j$.  Then by taking the
weak$^\ast$ limit, we have that \[f(z)\le \il{\partial K}
f(\zeta)\; d\om_K(z, \zeta).\]  Then $\om_K(z, \cdot)$
satisfies the sub-averaging inequality for every continuous
subharmonic function on $U_j$ and $\om_K(z, \cdot)$ is a
probability measure with support contained in $U_j$.  Thus
$\om_K(z, \cdot)$ must be in $\J^c_z(U_j)$, which is equal to
$\J_z(U_j)$ by Theorem~\ref{T:JensEQ}. Therefore $\om_K(z,
\cdot)\in \J_z(K)$.
\end{proof}

Following \cite[p. 16]{G78} a partial ordering on the set of Jensen measures is defined below.  The notation $\J(K)$ is used to stand for the union of all Jensen measures on $K$, that is
\[\J(K) = \bigcup_{z\in K} \J_z(K).\]

\bD
For $\mu,\nu\in\J(K)$ we say that $\mu \succeq \nu$ if for
every  $\phi\in S(K)$ we have
$\mu(\phi) \ge \nu(\phi)$.  Furthermore, a Jensen measure $\mu$ is maximal if there is no $\nu \succeq \mu$ with $\nu\ne\mu$
where $\nu\in \J(K)$.
\eD

We start with a simple observation.
\bL\label{L:indep}
If $\mu\in \J_{z_1}(K)$ and $\nu\in
\J_{z_2}(K)$ with $z_1\ne z_2$ then $\mu$ and $\nu$ are not
comparable.
\eL
\begin{proof}
To see this simply recall that the coordinate
functions $\pi_i$ are harmonic.  As $z_1\ne z_2$ they must
differ in at least one coordinate, say the $i^{th}$.  Assume
with out loss of generality that $\pi_i(z_1)>\pi_i(z_2)$.  Then
$\mu(\pi_i)>\nu(\pi_i)$.  However $-\pi_i$ is also harmonic and
so $\nu(-\pi_i)>\mu(-\pi_i)$.  Therefore $\mu$ and $\nu$ are
not comparable and if $\mu\succeq\nu$ then they have the common
barycenter.
\end{proof}

We will now show that the harmonic measure is maximal with respect to this ordering. The maximality of harmonic measure proved below is the Littlewood Subordination Principle (see \cite[Theorem 1.7]{D70}) when $K$ is the closed unit ball in the plane.
\bT\label{T:max}
For all $z\in K$, the measure $\omega_K(z,\cdot)$ is maximal in $\J(K)$.
\eT
\begin{proof}
By Lemma \ref{L:indep} it suffices to show that for any $z\in K$, $\omega_K(z, \cdot)$ is maximal in $\J_z(K)$.

Pick any $z_0\in K$.  Now we will show that $\omega_K(z_0,\cdot)$ majorizes every
measure $\mu\in\J_{z_0}(K)$. Consider a decreasing sequence of
regular domains $\{U_j\}$ converging to $K$.  Take any $\phi\in
S^c(K)$. By Proposition~\ref{P:SbhExt} we may find a sequence
$\phi_j\in S^c(U_j)$ increasing to $\phi$. Furthermore we
extend $\phi$ as $\tilde{\phi}\in C_0(\mathbb{R}^n)$ while
keeping $\tilde{\phi}\ge\phi_j$ for all $j$. Define harmonic
functions $\Phi_j$ on $U_j$ by
\[\Phi_j(x)=\int\limits_{\partial U_{j+1}}\phi_j(\zeta)\;d\omega_{U_{j+1}}(x,\zeta).\]
Therefore as $\phi_j$ is subharmonic, $\Phi_j\ge \phi_j$ on
$U_{j+1}$, so
\[\int\limits_{\partial U_{j+1}}\phi_j(\zeta)\;d\omega_{U_{j+1}}(z_0,\zeta)=\Phi_j(z_0)=
\mu(\Phi_j)\ge\mu(\phi_j).\] As $\tilde{\phi}\ge \phi_j$, we
have
\begin{equation}\int\limits_{\partial
U_{j+1}}\tilde{\phi}(\zeta)\,d\omega_{U_{j+1}}(z_0,\zeta)
\ge\mu(\phi_j),
\end{equation} for all $j$.  By taking
weak$^\ast$ limits, we have that
\[\lim_{j\rightarrow
\infty}\int\limits_{\partial
U_{j+1}}\tilde{\phi}(\zeta)\,d\omega_{U_{j+1}}(z_0,\zeta) =
\int\limits_{\partial K}\phi(\zeta)\;d\omega_K(z_0,\zeta).\] The Lebesgue Monotone Convergence
Theorem provides \[\lim_{j\rightarrow \infty}\mu(\phi_j) = \mu
(\phi).\] By taking the limit by $j$ of $(3)$ we see
\[\int\limits_{\partial K}\phi(\zeta)\;d\omega_K(z_0,\zeta)\ge \mu(\phi).\]
Therefore we have $\omega_K(z_0,\cdot)\succeq \mu$. If
any $\nu\in \J_{z_0}(K)$ has the property $\nu\succeq \mu$, by the antisymmetry property of partial orderings $\nu=\mu$. Thus the measure $\omega_K(z_0, \cdot)$ is maximal in $\J_{z_0}(K)$.
\end{proof}

The maximality of harmonic measures implies that they are
trivial at the points $z\in K$ such that $\J_z(K)=\{\delta_z\}$, which by Lemma \ref{L:Pfine} are precisely the fine boundary points.
\bC\label{C:h-triv}
The harmonic measure $\om_K(z_0,\cdot)=\delta_{z_0}$ if and
only if $\J_{z_0}(K)=\{\delta_{z_0}\}$.
\eC
\begin{proof}
Suppose $\om_K(z_0,\cdot)=\delta_{z_0}$.  Consider the function $\rho(z)=||z-z_0||^2\in S^c(K)$.  Then
for any $\mu\in \J_{z_0}$, by the maximality of
$\om_K(z_0,\cdot)$ we have
\[0=\rho(z_0)\le \mu(\rho)\le\int\limits_{\partial K} \rho(\zeta)\;d\om_K(z_0, \zeta)=\rho(z_0)=0.\]
As $\rho(z)>0$ for all $z\ne0$ and as $\mu$ is a probability
measure, we see that $\mu=\delta_{z_0}$.  Thus
$\J_{z_0}(K)=\{\delta_{z_0}\}$.

For the reverse implication we have already proved
Theorem~\ref{T:HJens} that $\om_K(z_0,\cdot)\in \J_{z_0}(K)$.
\end{proof}

\section{The Boundary}\label{S:Boundary}
In the book \cite{G78}, Gamelin introduces a version of Choquet theory for cones of functions on compact sets.  (Actually it applies to sets of functions which are slightly weaker than the cones we define.)

Following his guidance we consider a set $\mathcal{R}$ of functions mapping a compact set $K\subset \mathbb{R}^n$ to $ [-\infty, \infty)$ with the following properties:
\begin{enumerate}
\item[i.] $\mathcal{R}$ includes the constant functions,
\item[ii.] if $c\in \mathbb{R}^+$ and $f\in \mathcal{R}$ then $cf\in \mathcal{R}$,
\item[iii.] if  $f,g\in \mathcal{R}$ then $f+g\in \mathcal{R}$, and
\item[iv.] $\mathcal{R}$ separates the points of $K$.
\end{enumerate}

One then considers a set of $\mathcal{R}$-measures for $z\in K$ defined as the set of probability measures $\mu$ on $K$ such that
\[f(z) \le \mu(f)\]
for all $f\in \mathcal{R}$.

Naturally our model for $\mathcal{R}$ will be $S(K)$.  It then follows that when $\mathcal{R}=S(K)$ the $\mathcal{R}$-measures for $z\in K$ are precisely $\J_z(K)$.  We now state some classic results from \cite{G78} which we will need in the following sections.

One can define the Choquet boundary of $K$ with respect to $S(K)$ as
\[Ch_{S(K)}K = \{z\in K \colon \J_z(K) = \{\delta_z\}\}. \]
Many nice properties of the Choquet boundary are known.  In
particular, we will need the following characterization, see also,
for example, \cite[VI.4.1]{BH86} and \cite{H85}.

\bL\label{L:Pfine}  The Choquet boundary of $K$ with respect to $S(K)$
is the fine boundary of $K$, i.e.
\[Ch_{S(K)}K = \partial_f K.\]
\eL
\begin{proof}
Since the fine topology is strictly finer than the Euclidean
topology, any point in the interior of $K$ will also be in the
fine interior of $K$, and any point of $\mathbb{R}^n\setminus
K$  can be separated from $K$ by an Euclidean (therefore fine)
open set.  Therefore the fine boundary of $K$ is contained in
$\partial K$. The result follows immediately from \cite[Theorem 3.3]{P97} or \cite[Proposition 3.1]{BH86} which states that
$\mathcal{J}_z(K) = \{\delta_z\}$ if and only if the complement of $K$ is non-thin at $z$, that is $z$ is a fine boundary point of $K$.
\end{proof}

The set $\partial_fK$ is also called the stable boundary of
$K$.  In fact the lemma shows that $Ch_{S(K)}K$ is the finely
regular boundary of the fine interior of $K$.  For more details on finely regular boundary points and other related concepts,
see \cite[VII.5-7]{BH86} and \cite{H85}.

With this association, the result in \cite[p. 89]{B71} of Brelot about the stable boundary points of $K$ shows that $Ch_{S(K)}K$ is dense in $\partial K$.
\bT\label{T:DensPP} The fine boundary of $K$ (and therefore the Choquet boundary of $K$ with respect to $S(K)$) is dense in the topological boundary of $K$.
\eT

In general the fine boundary is not closed, as Example~\ref{E:SSC} of the Section \ref{S:DP} will show. So we cannot claim that it is the support of measures. Moreover, as Theorem \ref{T:DensPP} just showed the closure of $\rO_k$ is the boundary of $K$. In particular, it may coincide with $K$ for porous Swiss cheeses, see \cite[pg. 25-26]{G84}.

Recall that a measure $\mu \in \M(K)$ is concentrated on a set
$E$, if for every set $F\subset K\setminus E$, $\mu(F)=0$.  A
probability measure $\mu$ is concentrated on a set $E$ if and
only if $\mu(E)=1$.  From \cite[p. 19]{G78} we know that all maximal measures are concentrated on $Ch_{S(K)}K=\partial_fK$.  With this observation, the next corollary immediately follows from Theorem \ref{T:max} which stated that the harmonic measure is maximal.

\bC\label{C:conc}
For every $z$ in $K$, the harmonic measure with barycenter at
$z$ is concentrated on $\partial_fK$.
\eC

\section{The Dirichlet Problem on Compact Sets}\label{S:DP}

In the classical setting we know that any continuous function
in the boundary of a domain $D\subset \mathbb{R}^n$ extends
harmonically to $D$ and continuously to $\overline{D}$ if and only if every point of the boundary is regular.  For general compact sets in $\mathbb{R}^n$ we have the following result.

From this result it also follows that the swept-out point mass at $z$ onto $K$ is just $\om_K(z,\cdot)$.

\bT\label{T:regDP} If $K$ is a compact set in $\mathbb{R}^n$ then any function $\phi\in C(\partial_fK)$ extends to
a function in $H(K)$ if and only if the set
$\partial_fK$ is closed. Moreover, the solution is given by
\[\Phi(z) = \int\limits_{\partial_fK} \phi(\zeta)\; d\omega_K(z,\zeta) \qquad z\in K\] and $H(K)$
is isometrically isomorphic to $C(\partial_fK)$.
\eT
\begin{proof} Suppose that the set $\partial_fK$ is closed.
Consider a continuous function $\phi$ on $\partial_fK$.  Let
\[\Phi(z) = \int\limits_{\partial_fK} \phi(\zeta)\; d\omega_K(z,\zeta) \qquad z\in K.\]  As
$\partial_fK$ is closed, by Theorem~\ref{T:DensPP}, we have
$\partial_fK=\partial K$.  Also as $\omega_K(z, \cdot)=\delta_z$ for every
$z\in \partial_fK$, we see that $\Phi=\phi$ on $\partial_fK$.

Let $z_j$ be a sequence in $K$ converging to $z_0\in \partial_fK$.
As $z_0$ is in $\partial_fK = Ch_{S(K)}K$, so $\J_{z_0}(K)=\{\delta_{z_0}\}$.
Since (see \cite[p. 3]{G78}) $\J(K)$ is weak$^\ast$ compact, any sequence of measures $\mu_j\in \J_{z_j}(K)$ must converge weak$^\ast$ to $\delta_{z_0}$.  In particular, $\omega_{U_j}(z_j, \cdot)$ is
weak$^\ast$ converging to $\delta_{z_0}$.  Hence $\Phi(z_j)$ is
converging to $\Phi(z_0)=\phi(z_0)$, and $\Phi$ is continuous
at the boundary of $K$.

As $\partial_fK$ is closed, we have $\phi\in C(\partial_fK)=C(\partial K)$.
We extend $\phi$ continuously as $\tilde{\phi}\in
C_0(\mathbb{R}^n)$, and then define the harmonic functions
\[h_j(z) = \il{\partial U_j}\tilde{\phi}(\zeta) \;d\omega_{U_j}(z,\zeta).\]
As $\tilde{\phi}$ is continuous and $\omega_{U_j}(z,\cdot)$ converges
weak$^\ast$ to $\omega_K(z, \cdot)$,
\[\lim_{j \rightarrow \infty} h_j(z) =\lim_{j \rightarrow \infty} \il{\partial U_j}\tilde{\phi}(\zeta) \;d\omega_{U_j}(z,\zeta)= \il{\partial K}\phi(\zeta) \; d\omega_K(z,\zeta) = \Phi(z).\]

Therefore $\Phi$ is the pointwise limit of a sequence $\{h_j\}$ of functions harmonic in a neighborhood of $K$.  Furthermore we can take the extension $\tilde{\phi}$ of $\phi$ in such a way that the sequence $\{h_j\}$ is uniformly bounded.  It now easily follows that $\Phi$ is continuous on the interior of $K$.  Indeed, consider a
point $z$ in the interior of $K$.  Then there exists a ball $B$ centered at $z$ contained in the interior of $K$.  The $h_j$
are harmonic functions on $B$ and converging pointwise to
$\Phi$.  Thus $\Phi$ is continuous on $B$ by the Harnack
principle, and so $\Phi$ is continuous on $K$.  Therefore we
have a continuous function $\Phi$ with representation
\[\Phi(z) = \int\limits_{\partial K} \phi(\zeta)\; d\omega_K(z,\zeta)  \qquad z\in K.\]

Since $\Phi$ is continuous on $K$ by \cite{P97} to check that $\Phi\in H(K)$ all that remains is to show that $\Phi$ is averaging with respect to Jensen measures, i.e. the equivalence of the external definition of $H(K)$ and the definition by Jensen measures. So we need to see that  $\Phi(z)=\mu_z(\Phi)$ for every $\mu_z \in \J_z(K)$ and for every $z\in K$.  As $h_j$ is harmonic on $U_j$, $h_j(z)=\mu_z(h_j)$.  However by the Lebesgue Dominated Convergence Theorem
\[\mu_z(\Phi)=\lim_{j\rightarrow\infty}\mu_z(h_j)=\lim_{j\rightarrow\infty}h_j(z)=\Phi(z).\]
Thus $\Phi\in H(K)$.

For the converse, suppose $\partial_fK$ is not closed.  Then there is a point $z_0\in\partial K \setminus \partial_fK$.  Since $z_0$ is not in $\partial_fK$, by Corollary \ref{C:h-triv}, $\om_K(z_0,\cdot)$ is not trivial.  Therefore we can find a set $E\subset \partial K$ such that $\om_K(z_0, E)>0$ with $E$ in the complement of $B(z_0, r)$ for some $r>0$. Consider a continuous function $f$ on $\partial K$ such that
$f=1$ on $\bd K$ outside $B(z_0, r)$ is $1$ and $f=0$ on
$B(z_0, r/2)\cap \partial K$. Then
\[\int\limits_{\partial K} f(\zeta)\; d\om_K(z_0, \zeta)>\om_K(z_0, E) \qquad z\in K.\]
However $f(z_0)=0$.  Thus there can be no function in $H(K)$
which agrees with $f$ on the boundary of $K$.
 \end{proof}

\begin{Example}\label{E:SSC}
The following set provides a simple example of a compact set $K\subset \mathbb{R}^n$, $n\ge 3$, in which the fine boundary is not closed.  The set $K$ is obtained from the closed unit ball $\overline{B}\subset \mathbb{R}^n$ by deleting a sequence $\{B(z_n, r_n)\}_{n=1}^\infty$ of open balls whose centers and radii tend to zero.  We take the centers to be $z_n=(2^{-n}, 0, \ldots, 0 )\in \mathbb{R}^n$ and the radii $0<r_n<2^{-n-2}$.  This example is analogous to the ``road runner'' example of Gamelin \cite[Figure 2, pg 52]{G84} and the Lebesgue spine \cite[pg 187]{AG01}.
\end{Example}

By Theorem \ref{T:regDP} one can not expect a continuous
solution for the Dirichlet problem on an arbitrary compact set
even with continuous boundary data. Therefore at this point we
consider the following broader class of solutions with weaker
continuity requirement.

\bD
Let $fH^c(K)$ denote the class of finely continuous functions on $K$ which are finely harmonic on the fine interior of $K$ and continuous and
bounded on $\partial_fK$.
\eD

We have seen (the definition via Jensen measures) that $H(K)$ consists of the functions in $C(K)$ satisfying the averaging property with respect to $\J(K)$ and by the interior definition of $H(K)$ can also be seen as the $C(K)$ functions
which are finely harmonic on the fine interior of $K$.
Therefore in the definition of $fH^c(K)$ we have maintained the finely harmonic requirement while requiring continuity only on the boundary $\partial_fK$ (to match the boundary data).  In fact Theorem~\ref{T:qdp} below shows that the functions in $fH^c(K)$ also satisfy the averaging property with respect to $\J(K)$.

Theorem~\ref{T:qdp} will show that the Dirichlet problem on
compact sets $K\subset \mathbb{R}^n$ is solvable in the class
of functions $fH^c(K)$ for boundary data that is continuous and bounded on $\partial_fK$.  The functions which are continuous and bounded on $\partial_fK$ will be denoted $C_b(\partial_fK)$.  For this we will need the following \cite[Theorem 11.9]{F72} of
Fuglede.
\bT\label{T:Fugl} The pointwise limit of a pointwise convergent sequence of finely harmonic functions $u_m$ in $U$, a finely open subset of $\mathbb{R}^n$, is finely harmonic provided that $\sup_m |u_m|$ is finely locally bounded in $U$.
\eT

\bT\label{T:qdp}
For every $\phi\in C_b(\partial_fK)$, i.e. continuous and bounded on $\partial_fK$, there is a unique $h_\phi \in fH^c(K)$ equal to $\phi$ on $\partial_fK$.  Moreover, $h_\phi$ satisfies the averaging property for $\J(K)$ and in particular
\[h_\phi(x)=\int\limits_{\partial_fK}\phi(\zeta) \; d\omega_K(x,\zeta), \qquad x\in K.\]
\eT
\begin{proof}
Let $\phi\in C_b(\partial_fK)$ and for $x\in\overline{\partial_fK}$ define
\[\tilde{\phi}(x)=\limsup_{y\rightarrow x,
\; y\in \partial_fK} \phi(y).\]  Since $\phi$ is
continuous on $\partial_fK$, if $x\in \partial_fK$ then
$\tilde{\phi}(x)=\phi(x)$.  Furthermore, $\tilde{\phi}$ is
upper semicontinuous, and as such we may find a decreasing
sequence of functions $\{\phi_k\}$ which are continuous on
$\overline{\partial_fK}$ and converge pointwise to $\tilde{\phi}$.
Then we extend the $\phi_k$ to $C_0(\mathbb{R}^n)$ as
$\hat{\phi}_k$. By taking $\tilde{\phi}_k= \min\{\hat{\phi}_1,
\hat{\phi}_2, \cdots ,\hat{\phi}_k\}$ we can make the
extensions be decreasing.  Consider a decreasing sequence of
regular domains $U_j$ converging to $K$.  Let $u_{j,\,k}$ be
the solution of the Dirichlet problem on $U_j$ for
$\tilde{\phi}_k$.  As the measures $\omega_{U_j}(x,\cdot)$
weak$^\ast$ converge to $\omega_K(x, \cdot)$, we have that
$\lim_{j}u_{j,\,k}=\int \tilde{\phi}_k\;d\omega_K
\mathrel{\mathop:}= u_k$. As the $ \tilde{\phi}_k$ are
decreasing, $u_k$ must also be decreasing.  Indeed, we will let $h_\phi=\lim u_k$.

Take any $\mu\in \J(K)$.  Then $\mu\in \J_{z_0}(U_j)$ for all $j$ and some $z_0\in K$.  As $u_{j,\,k}$ is harmonic, we have $\mu(u_{j,\, k})=u_{j,\,k}(z_0)$.  However by the Lebesgue Dominated Convergence Theorem we have $\lim_j \mu(u_{j,\, k}) = \mu(u_k)$, and so $\mu(u_k) = u_k(z_0)$.  Since the sequence $\{u_k\}$ is decreasing pointwise to $h_\phi$ we have that $\mu(h_\phi)=h_\phi(z_0)$ by the Lebesgue Monotone Convergence Theorem.  Thus $h_\phi$ satisfies the averaging property on $\J(K)$.  As $\om_K(z, \cdot)\in \J(K)$ for all $z\in K$ we see that
\[h_\phi(z) = \int\limits_{\partial_fK} h_\phi(\zeta) \; \om_K(z, \zeta).\]
We will now show that $h_\phi = \phi$ on $\partial_fK$.  For any $x\in \rO_k$, we know $\omega_K(x,\cdot)=\delta_x$, and
\[u_k(x)=\lim_{j\rightarrow \infty}u_{j,\,k}(x)=\int \tilde{\phi}_k(\zeta)\;d\omega_K(x,\zeta) = \tilde{\phi}_k(x).\]
Thus $u_k(x)=\tilde{\phi}_k(x)$ for all $x\in \partial_fK$, and so
\[h_\phi(x)=\lim_{k\rightarrow \infty} u_k(x)=\lim_{k\rightarrow \infty} \tilde{\phi}_k(x) =\phi(x),\]
 for all $x\in\partial_fK$.

To see that $h_\phi$ is finely harmonic we use
Theorem~\ref{T:Fugl}.  Observe that $u_k$ is the pointwise
limits of the harmonic (and therefore finely harmonic)
functions $u_{j,\, k}$, and the solution $h_\phi$ is the
pointwise limit of $u_k$. From the construction of these
functions it is clear that they are bounded.
\end{proof}

\bC  The set $C_b(\partial_fK)$ is isometrically isomorphic to $fH^c(K)$.
\eC
\begin{proof}
The previous theorem establishes the homomorphism taking
$C_b(\partial_fK)$ to $fH^c(K)$.  Observe that $h|_{\partial_fK}\in
C_b(\partial_fK)$  for every $h\in fH^c(K)$.  The uniqueness of the
solution shows that $h|_{\partial_fK}$ extends as $h$.  Furthermore,
the isometry follows directly from the integral representation
in the previous theorem.
\end{proof}

\end{document}